# WheelCon: A wheel control-based gaming platform for studying human sensorimotor control

Quanying Liu[1,2,3*], Yorie Nakahira[1], Ahkeel Mohideen[1], Adam Dai[1],

Sung Hoon Choi[1], Angelina Pan[1], Dimitar M. Ho[1], John C. Doyle[1*]

[1] Division of Engineering and Applied Science, California Institute of Technology, Pasadena, CA, USA

[2] Neuroscience Center, Huntington Medical Research Institutes, Pasadena, CA, USA

[3] Research Center for Motor Control and Neuroplasticity, KU Leuven, Leuven, Belgium

*Corresponding author:

Quanying Liu, quanying.liu@caltech.edu

John C. Doyle, doyle@caltech.edu

**Abstract**

Feedback control theory has been extensively implemented to theoretically model human sensorimotor control. However, experimental platforms capable of manipulating important components of multiple feedback loops lack development. This paper describes WheelCon -- an open source platform aimed at resolving such insufficiencies. Using only a computer, standard display, and inexpensive gaming steering wheel equipped with a force feedback motor, WheelCon safely simulates the canonical sensorimotor task of riding a mountain bike down a steep, twisting, bumpy trail. The platform provides flexibility, as will be demonstrated in the demos provided, so that researchers may manipulate the disturbances, delay, and quantization (data rate) in the layered feedback loops, including a high-level advanced plan layer and a low-level delayed reflex layer. In this paper, we illustrate WheelCon's graphical user interface (GUI), the input and output of existing demos, and how to design new games. In addition, we present the basic feedback model and the experimental results from our demo games, which align well with the model's prediction. The WheelCon platform can be downloaded at https://github.com/Doyle-Lab/WheelCon. In short, the platform is featured as cheap, simple to use, and flexible to program for effective sensorimotor neuroscience research and control engineering education.

# 1. Introduction

Human sensorimotor control system is extremely robust (Franklin and Wolpert, 2011), although the sensing is distributed, sparse, quantized, noisy and delayed (Bays and Wolpert, 2007; Desmurget and Grafton, 2000; Sanger and Merzenich, 2000); the computing in the central nervous system is slow (Mohler and Okada, 1977; Muller, et al., 2018; Zhang, et al., 2018); and the muscle actuation fatigues and saturates (Blinks, et al., 1978). Many theoretical models have been proposed to explain the complicated human sensorimotor control process (Gallivan, et al., 2018; Sanger, 2018; Todorov, 2004; Todorov and Jordan, 2002). For example, feedback control theory predicts the optimal control policy (Todorov and Jordan, 2002), and Bayesian theory models sensorimotor learning (Kording and Wolpert, 2004). In contrast to the abundance of theoretical models, experimental platforms capable of manipulating important components of multiple feedback loops lack development. This is in part due to the fact that designing a platform to bridge and test these aspects of sensorimotor control requires a diverse range of expertise, extending from motor control theory, signal processing and interaction, all the way to computer graphics and programming. Researchers often develop their own custom hardware/software systems to characterize human sensorimotor control performance, which can limit the ability to compare/contrast and integrate datasets across research groups. Development of an easy-to-use and validated system could broaden the quantitative characterization of sensorimotor control.

In this paper, we present the WheelCon platform, a novel, free and open source platform to design video games for a virtual environment that non-invasively simulates riding a mountain bike down a steep, twisting, bumpy trail.

WheelCon contains the highly demanded basic components which present in each theory: delay, quantization, noise, disturbance, and multiple feedback loops. It is a potential tool for studying the following diverse questions in human sensorimotor control:

- How the human sensorimotor system deals with the delay and quantization in neural signaling, which are fundamentally constrained by the limited resources (such as the space and metabolic costs) in our brain (Nakahira, et al., 2015);
- How humans handle the unpredictable, external disturbances in sensorimotor control (Wolpert and Miall, 1996);
- How the hierarchical control loops layered and integrated in human sensorimotor system (Battaglia-Mayer, et al., 2003; Scott, 2004);
- The consequence of the delay and quantization in human visual feedback (Saunders and Knill, 2003) and reflex feedback (Insperger, et al., 2013) in sensorimotor control;
- The optimal policy and strategy for sensorimotor learning under delay and quantization (Kording and Wolpert, 2004).

WheelCon integrates with a steering wheel and can simulate game conditions which manipulate the variables in these questions, such as signaling delay, quantization, noise and disturbance, while recording the dynamic control policy and system errors. It also allows researchers to study the layered architecture in sensorimotor control. In the example of riding a mountain bike, two control layers are involved in this task: the high-layer plan and the low-layer reflex. For visible disturbances (i.e. the trail), we plan before the disturbance arrives. For disturbances unknown in advance (i.e. small bumps), our control relies on delayed reflexes. Feedback control theory proposes that effective layered architectures can integrate the higher layers' goals, plans, decisions with the lower layers' sensing, reflex, and action (Nakahira, et al., 2015). WheelCon provides experimental tools to induce distinctive disturbances in the plan and reflex layers separately for testing such a layered architecture (Figure 1).

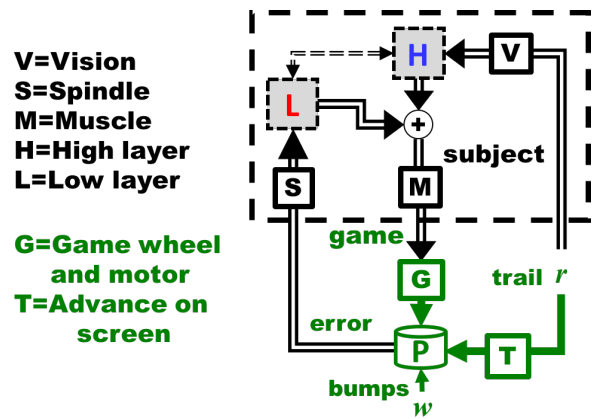

Figure 1 - Basic block diagram for experimental platform with subject and gaming wheel with motor. Each box is a component that communicates or computes and has potentially both delay and quantization, including within the game in G. The advance warning T is also implemented on a computer screen with vision.

The platform included the visual input, sensorimotor control and visual feedback. Using these components, we have implemented several demo video games for the users to quantitatively measure human sensorimotor control performance in multiple scenarios, which connect the tasks with real and virtual environment. This paper is organized as follows: Section 2 describes the features of WheelCon; Section 3 presents two simple examples of using basic feedback control theory to model the delay and quantization in a single control loop; Section 4 presents five demo games and the corresponding experimental results. Finally, section 5 summarizes the significance of WheelCon from the neuroscience perspective and the control education perspective.

## 2. Platform Description

The WheelCon platform can be downloaded at https://github.com/Doyle-Lab/WheelCon. It includes several demo games, code for low-level development, and the source code for high-level development.

**The Hardware**

WheelCon requires a gaming wheel as the controller, a monitor to present the visual input and output, and a computer to decode the controller's commands and to synchronize with the game.

We have tested the platform with a steering wheel manufactured by both Logitech and Thrustmaster. Different wheels have different properties like torque capability, so the game's input file must be tuned to match desired force feedback effects.

Any monitor should be capable of properly displaying the visual input, although a minimum screen size of 15 inches and a resolution of 1920 x 1080 is recommended. The game must be run on a Windows 10 computer with preferably high-performance capabilities. This includes a recommended i7 Intel processor or equivalent AMD chip and a minimum of 8 GB of RAM.

**The Control Variables in WheelCon**

WheelCon allows users to manipulate the multiple control components, such as the disturbances, delay and quantization, of the layered feedback loops, which includes an advanced plan layer and a delayed reflex layer. The variables which can be manipulated in control loops are shown in Table 1.

These variables are set in a .txt input file, making it easy for users to design new games without editing the source code. Upon startup, the game loads the variables' values into Animation Curves for the software to track. This means that at every time, every experimental variable must be associated with a value. The resolution for the time is at minimum 10ms (0.01 seconds). Therefore, in the input file, all experimental variables must be valued at every time stamp. Each line of the text file should represent a single time stamp with appropriate variable values. The format of the line depends on the version of the game being played. For the mountain task the format is: time, horizontal position of the trail disturbance, bump disturbance, action data rate, action delay, vision delay, vision quantization (e.g. 0.01,6,10,-1,30,0.2).

Table 1 – The variables in control loops which can be manipulated

| notation | variable | unit | constrains |
|---|---|---|---|
| w(t) | Bump disturbance | n Newton | $0 \leq w(t) \leq 100$ |
| r(t) | Trail disturbance | 100 pixels | $0 \leq r(t) \leq 100$; $r(t) \perp w(t)$ |
| $T_{vis}$ | Vision advance warning/delay | second | $-1 \leq T_{vis} < 1$ |
| $T_{act}$ | Action delay | second | $0 \leq T_{vis} < \infty$ |
| $R_{vis}$ | External data rate in vision input | bit | $1 \leq R_{vis} \leq 10$ |
| $R_{act}$ | External data rate in the wheel output | bit | $1 \leq R_{act} \leq 10$ |
| $Q_{vis}$ | External quantizer on the vision input | 1 | $Q_{vis} = 2^{R_{vis}}$ |
| $Q_{act}$ | External quantizer on the wheel output | 1 | $Q_{act} = 2^{R_{vis}}$ |

**The Output from WheelCon**

Besides the control variables in the input file, the error and control policy will also be exported for each sampling time. See Table 2 for the detailed description.

Table 2 – The variables in control loops which can be manipulated

| notation | variable | unit |
|---|---|---|
| t | Time | second |
| x(t) | Error dynamics, measured by the distance between the player's position and the desired position | 100 pixels |
| u(t) | Control policy, measured by the deviated angle of the steering wheel | degree |

**GUI of WheelCon**

Download WheelCon, and follow the installation in manual. Start the platform by double clicking on 'WheelCon.exe'.

WheelCon has a main menu as seen in Figure 2(a). From here, the sensitivity of the wheel for both games can be altered and whatever game to played can be selected.

Clicking on a game button brings up either Figure 2(b) or Figure 2(c). These submenus represent the game setup for whichever task; here, subjects can enter their name for the output file and researchers select the correct game input file.

Fitt's Law Task: The video game's interface for the Fitt's Law Task is shown in Figure 2(d). The Fitt's Law Task mimics a traditional reaching task; the player turns the steering wheel left or right, in order to reach a gray zone. The player should stay in the zone until the zone jumped to another area. The movement time is measured as the time between the moment the grey zone changes and the moment the green bar enters the grey zone.

Mountain Bike Task: The mountain bike task safely simulates riding a mountain bike down a steep, twisted, bumpy trail using a standard display and a gaming steering wheel. The virtual trail scrolls down a PC screen which can vary in speed, turns, and visual look ahead (and thus advanced warning or delay). Subjects can see the trail and turn the wheel to track it with minimum error, while an internal motor can torque the wheel to mimic invisible bumps in the trail. The video game interface for mountain bike task is shown in Figure 2(e). The grey line is the desired path, the green line is the player's current position. The aim of the task is to control the green line so that it tracks the grey line with minimal error, where error is the distance between the green line and grey line at the present time.

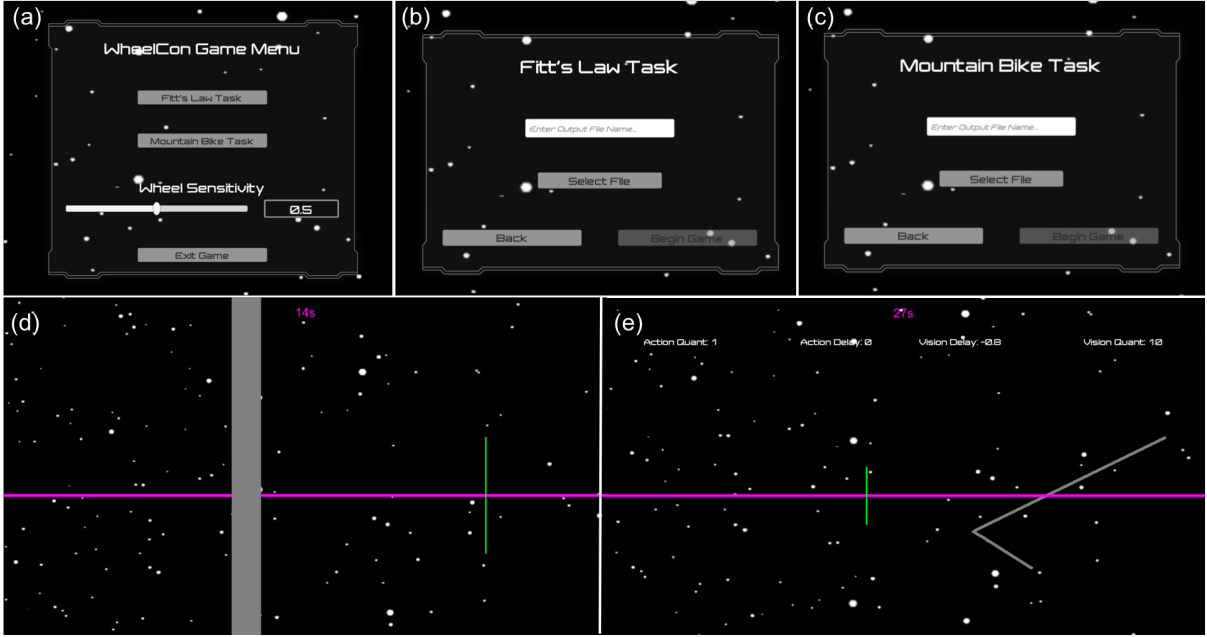

Figure 2 – The user-graphic interface for WheelCon: (a) the main menu; (b) the Fitt's Law Task menu; (c) the Mountain Bike Task menu; (d) the video game interface for Fitt's Law Task; (e) the video game interface for Mountain Bike Task.

## 3. Basic Feedback Control Model

We show a simplified feedback control model shown in Figure 4. The system dynamics is given by

$$x(t+1) = x(t) + w(t) + u(t) + r(t)$$

where x(t) is the error at time t, r(t) is the trail disturbance, w(t) is the bump disturbance, and u(t) is the control action.

**Modelling Action Delay in Trail Disturbance**

When there is a delay T in the action, and a trail disturbance r(t), we model the control action by

$$u(t + T) = \kappa\big(x(0{:}t), r(0{:}t), u(0{:}t + T - 1)\big). \quad (1)$$

The game starts with zero initial condition, i.e., $x(0) = 0$. The controller $\kappa$ generates the control command $u(t)$ using the full information on the histories of state, disturbance, and control input. Here, the net delay T is composed of the internal delays in the human sensorimotor feedback and the delays externally added. The control command is executed with delay T ≥ 0.

Sensorimotor control in the risk-aware setting motivates the use of L1 optimal control, and as such, our goal is to verify the following robust control problem

$$\inf_{\kappa} \sup_{\|r\|_\infty \leq 1} \|x\|_\infty . \quad (2)$$

This problem admits a simple and intuitive solution. The optimal cost is given by

$$\inf_{\kappa} \sup_{\|r\|_\infty \leq 1} \|x\|_\infty = T . \quad (3)$$

This optimal cost is achieved by the worst-case control policy $u(t + T) = -r(t)$, which yields

$$\inf_{\kappa} \sup_{\|r\|_\infty \leq 1} \|u\|_\infty = 1 . \quad (4)$$

**Modelling Action Quantization in Trail Disturbance**

When the data rate, R, in the control loop is limited, the control action is generated by the following feedback loop with communication constrains,

$$u(t) = Q(\kappa_t(x(0:t), u(0:t-1))) , \quad (5)$$

where $\kappa_t : (\mathbb{R}^{t+1}, \mathbb{R}^t) \to \mathbb{R}$ is a controller, and $Q: \mathbb{R} \to S$ is a quantizer with data rate R ≥ 1, i.e. S is a finite set of cardinality $2^R$. The disturbance r(t) is infinity-norm bound and without loss of generality, $\|r\|_\infty \leq 1$.

The worst-case state deviation is lower-bounded by

$$\sup_{\|r\|_\infty \leq 1} \|x\|_\infty \geq \frac{1}{2^R - 1} , \quad (6)$$

and the minimum control effort is given by

$$\sup_{\|r\|_\infty \leq 1} \|u\|_\infty \geq \left(1 + \frac{1}{2^R - 1}\right)\left(1 - \frac{1}{2^R}\right) . \quad (7)$$

**Measures of errors**

To quantify the performance, we measured the $L_\infty$- norm, L1-norm and L2-norm of error. The infinity norm, $L_\infty$- norm, is defined as the maximum of the absolute errors, where

$$\|x\|_\infty = max(|x(1)|, |x(2)|, \ldots, |x(n)|) \ . \ (8)$$

L1-norm is known as absolute errors. Here we define L1-norm as the mean of the sum of the absolute error as following

$$\|x\|_1 = \frac{1}{n}\sum_{t=1}^{n}|x(t)|. \ (9)$$

L2-norm is the square root of the sum of the squares of the absolute errors. Here we define the L2-norm error as following

$$\|x\|_2 = \sqrt{\frac{1}{n}\sum_{t=1}^{n}|x(t)|^2}. \ (10)$$

## 4. Experimental results from WheelCon platform

We designed specific game trials for WheelCon to show the potential applications. For example, the platform can be used to examine how well a player can drive along a given trail in the virtual environment with control/communication constraints, such as limited visual advanced warning while being exposed to force disturbances from the internal motor. The research protocol was approved by the California Institute of Technology IRB and the subject provided informed consent prior to any procedures being performed. To show the capability of WheelCon, we presented the tests and results from one subject below.

**Game 1: Visual Advanced Warning or Delay**

Game 1 of WheelCon evaluates how the length of the look-ahead window (advanced warning / delay) affects sensorimotor control performance without being exposed to additional disturbances.

Game 1 lasts for 400 seconds and consists of one continuous "Trail" which reduces the amount of look-ahead every 30 seconds. The game begins with 1 s of advanced warning, then decreases to 0.75 s, and then to 0.5 s. From there, the game decreases the look ahead by 0.1 s until a minimum of -0.5 s is reached. Positive delay, or negative advanced warning, means only the trail behind the player is visible.

An evolution of the absolute error of the player as the game progresses is depicted in Figure 3(a). The plot displays only the middle 20 seconds of each of the 30 second intervals to neglect effects of the player adjusting to the new look-ahead window. The progression of the error in the blocks of constant look-ahead demonstrate the intuitive relationship, that the player loses performance with a decrease of advanced warning. To quantify that effect in more detail, we evaluate $L_1$-/$L_2$-/$L_\infty$- norm for the error dynamics for every 20 second group corresponding to a delay level. Summarizing these calculations in a plot gives us Figure 3(b) which demonstrates how the players error-norm does not change until the advanced warning reaches 0.3 seconds and then increases in an approximately linear fashion.

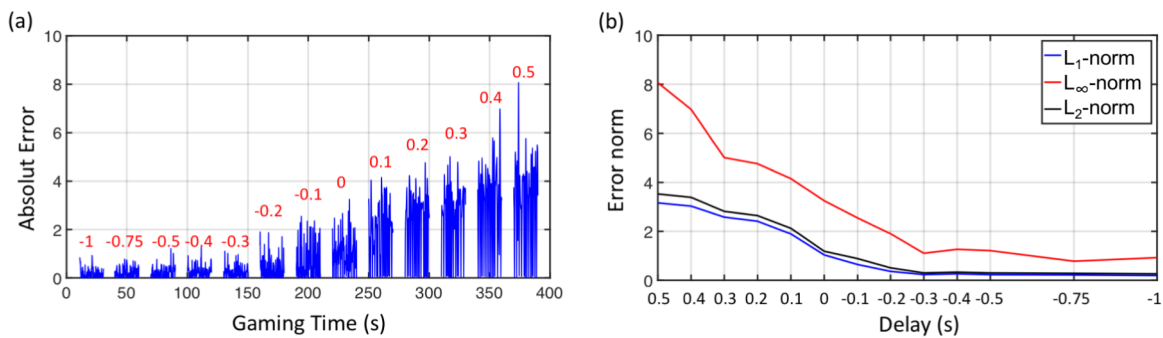

Figure 3 – Adding delay in vision input during the mountain-bike task. (a) The absolute error changes with the gaming time. The corresponding delay/advanced warning is shown with the red number. (b) The $L_1$-/$L_2$-/$L_\infty$- norm of error decreases with the decreasing delay. The negative value for delay means advanced warning.

**Game 2: Delay in Action Output**

Unlike Game 1's external visual delay, Game 2 adds specific internal delay to the action output; in other words, the current control policy u(t) works at u(t+$T_{act}$) where $T_{act}$ is the external delay in action. Game 2 lasts for 180 seconds. Adjusting $T_{act}$ every 30 seconds, $T_{act}$ starts at 0 s, and increments by 0.15 s until it reaches 0.75 s.

The effects of delay in action are shown in Figure 4. Similar to the vision delay, the error increases linearly with the delay, which is well in line with the prediction from theory in Eq(3).

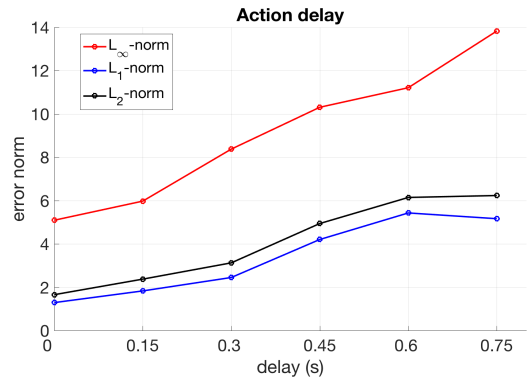

Figure 4 – Adding delay in action output. The L1-/L2-/L∞- norm of error increases with the increasing delay.

**Game 3 and 4: Quantization in Vision input and Action Output**

Game 3 and Game 4 study the effects of quantization in vision input and action output respectively. Each game is 210 seconds long and the quantization changes every 30s with the data rate increasing from 1 to 7 bits. For example, when the $R_{vis}$ is 1 in Game 3, the desired position (gray line in the gaming GUI) is presented either in the center left or the center right of the screen. When $R_{vis}$ is n, the desired position can be presented in $2^n$ possible locations in the screen. For Game 4, when $R_{act}$ is 1, the player is either going left or right with one speed. When $R_{act}$ is n, the player can steer the wheel to go left or right with $2^{n-1}$ speeds.

The effects of quantization (limited data rate) in the vision and action are shown in Figure 5. Consist with the theory's prediction in Eq(6), sensorimotor control performance improves with higher data rates, and reaches the optimal control performance when R is around 5.

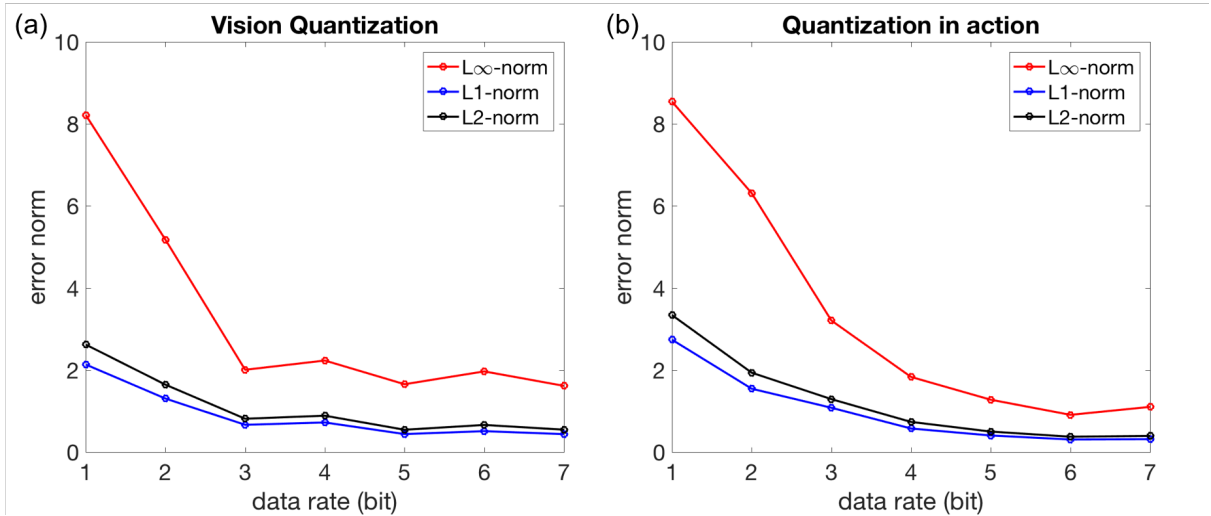

Figure 5 – Quantization in vision input (a) and action output (b). The L1-/L2-/L∞- norm of error are shown in the blue, black and red line respectively.

**Game 5: Bump and trail disturbance**

Game 5 is designed to test the effects of bump and trail disturbances on human sensorimotor control. Game 5 consists of three scenarios:

a) "Bumps", tracking a constant trail subject despite torque disturbances on the wheel that mimic hitting bumps when riding a mountain bike;

b) "Trail", tracking a moving trail with random turns but without bumps;

c) "Trail with Bumps", tracking a moving trail with random turns and bumps.

Each scenario lasts for 60 seconds in the order (Bumps, Trail, Trail with Bumps) with 5 seconds rest preceding each scenario. Furthermore, the disturbances and the trail during the isolated phases are duplicated in the combined "Trails with Bumps" phase, so that a proper performance comparison can be drawn between the separate tasks and the one where the player must multiplex. During the entire game, there is 1s of advanced warning in vision input, no delay in action output, and a 10-bit data rate for both vision and action.

As the disturbance, we use a random, binary signal, whose amplitude is the maximum possible torque the motor of the steering wheel can exert. Every 100 ms, the torque switches between max positive and negative (100 or -100 for the Logitech G27 wheel).

A similar random binary switching controls the trail derivative. More specifically, the trail travels at a constant speed but randomly switches its direction such that it always stays in the screen range comfortably visible to the player. We adjusted the velocity of the trail on the screen such that the required steering wheel turning rate is approximately 75° / second.

Figure 6 illustrates 5 second snapshots of the error dynamics for each scenario during the game.

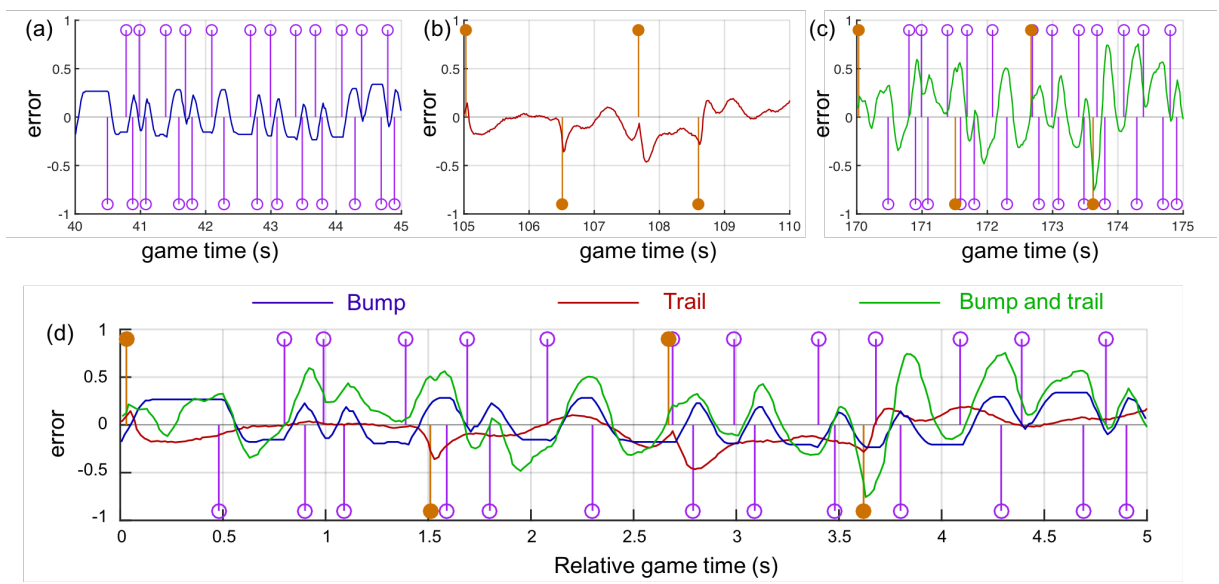

Figure 6 – Effects of bump and trail disturbance on human sensorimotor control. (a) the error dynamics induced by bump disturbance; (b) the error dynamics induced by trail disturbance; (c) the error dynamics induced by bump and trail disturbance; (d) the overlayed error in bumps (blue), trails (red), Trails with bumps (green). The purple-empty and orange-filled stem plots indicate the timing and direction of the bump disturbances and trail disturbance, respectively. Note that both the wheel forces and the trail rates are square waves, and the stems indicate where these square waves switch (i.e. derivatives of the forces and rates).

## 5. Discussions

In this paper, we have presented a free, open-source gaming platform, WheelCon, for studying the effects of delay, quantization, disturbance, and layered feedback loops in human sensorimotor control. Our platform provides the necessary tools to non-invasively manipulate external delays and limit data rate in both vision inputs and action output. We have showed the hardware, software, GUI. The settings of a single sensorimotor control loop with delay and quantization has been implemented, which allows us to measure the effects of delay, quantization and disturbance in sensorimotor control. The experimental results are well in line with the prediction from the feedback control theory. The executable file, demo inputs, source code, and detailed manual of WheelCon can be downloaded at https://github.com/Doyle-Lab/WheelCon.

**Neuroscience Perspective**

The demo games in the WheelCon platform showed how to use WheelCon to study the effects of delay (Game 1 and 2), quantization (Game 3 and 4) and disturbance in advanced plan layer and delayed reflex layer (Game 5) in human sensorimotor motor control. However, these examples are only the beginning of a new epoch. We serve much broader interests: to study the speed-accuracy tradeoffs (MacKay, 1982; Osman, et al., 2000), the layered architectures (Doyle and Csete, 2011; Nakahira, et al., 2015) and the learning strategy (Haruno, et al., 2001; Kording and Wolpert, 2004) in human sensorimotor control combined with mathematical models from information theory, feedback control theory, reinforcement learning or Bayesian inference. For instance, the data presented in this paper is only recorded with optimal control and minimal training, we did not analyze the data during learning. Some specific games can be designed for studying the sensorimotor learning under the inevitable delay, limited data rate, large initial error and unpredictable disturbance (Sanger, 2004; Schneider and Detweiler, 1988; Shmuelof, et al., 2012).

Along this line, we would like to propose some potential research questions that could be studied with our experimental platform in the future:

- Is the optimal control policy obtained from control theory shown in Eq(4) and Eq(7) applied in real practice? In this paper, we only examined system performance. It is valuable to investigate the control policy in human sensorimotor control.

- How do the tradeoffs in speed / accuracy / saturation / energy cost in muscles benefit from the diversity of neurons? For example, we can set the speed of actuator (wheel) as a uniform speed to mimic the uniform muscle; or allow multiple speeds from the wheel to mimic the diverse types of muscles. The saturation of muscles can be simulated as a maximum speed generated by steering wheel.
- How do the speed-accuracy tradeoffs in the high-layer plan/decision making support human/animals' decision strategies across complex environment under uncertainty, limited information, and risks?
- How does the human sensorimotor system tolerate the noise in control loops? It can be tested by adding noise to actuator or vision (the trajectory of trail shown in the screen).
- What are the effects of learning/adaptation in different control layers? Do the tradeoffs exist between fast learning and fast forgetting, between efficiency and plasticity?

We provide a cheap, easy to use and flexible to program platform, WheelCon that bridges the gap between theoretical and experimental studies on neuroscience. Specifically, it can be used for examining the effects of delay, quantization, disturbance, potentially speed-accuracy tradeoffs. It can also be applied for studying decision making, and multiplexing ability across different control layers in human sensorimotor control. Moreover, WheelCon is compatible with non-invasive neural recordings, such as electroencephalography (EEG), to measure the neural response during sensorimotor control (Birbaumer, 2006; Liu, et al., 2017; Nicolelis, 2003), and the non-invasive brain stimulation techniques, such as Transcranial Electrical Stimulation (tES) and Transcranial Magnetic Stimulation (TMS), to manipulate the neural activity (Hallett, 2000; Madhavan, et al., 2011).

**Control engineering perspective**

WheelCon allows users to design specific disturbances and noise in control loops, to add delay and quantization in both visual sensing and action output. It will also provide potential measures to quantify system performance (i.e. $L_1$-norm, $L_2$-norm and $L_\infty$- norm), input and actuator output. These are essential elements in a control system, and WheelCon's capability of analyzing them suggests its potential use in control education.

Robust control theory is a powerful tool in the investigation on the effect of noise, disturbances, and other uncertainties in system performance (Leong and Doyle, 2016). However, the impact of such theory is limited by its technical accessibility. Given the platform, we can easily study and compare the settings with delay and quantization with those of delay or quantization. This clear separation of constraints in the feedback loop can help explain the basics of control in uncertain dynamical systems, and allow us to demonstrate how the plant instability, actuator saturation, and unstable zero dynamics impact our sensorimotor system.

In summary, WheelCon is a versatile platform to design diverse games for studying the effects of disturbance, delay and quantization in the layered control loops in human sensorimotor control. We believe it will address crucial gaps in both neuroscience and control education by demonstrating the robustness, reliability, and efficiency of human sensorimotor control.

**Acknowledgement**

We thank CIT Endowment & National Science Foundation (to J.C.D) and Boswell postdoctoral fellowship and FWO postdoctoral fellowship (12P6719N LV to Q.L) for the supports.

**Conflict of Interest**

The authors declare that they have no conflicts of interest.